\documentclass[12pt]{amsart}

\usepackage[text={425pt,650pt},centering]{geometry}
\usepackage{color}
\usepackage{esint,amssymb}
\usepackage{graphicx}
\usepackage{tikz}
\usetikzlibrary{calc} 
\usepackage{slashed}
\usepackage{mathtools}
\usepackage{dsfont}

\newtheorem{theorem}{Theorem}
\newtheorem{proposition}[theorem]{Proposition}

\theoremstyle{definition}

\newcommand{\cref}[1]{Corollary~\ref{c.#1}}

\numberwithin{equation}{section}
\numberwithin{theorem}{section}

\newcommand{\Z}{\mathbb{Z}}
\newcommand{\N}{\mathbb{N}}

\newcommand{\R}{\mathbb{R}}

\newcommand{\E}{\mathbb{E}}
\renewcommand{\P}{\mathbb{P}}
\newcommand{\F}{\mathcal{F}}
\newcommand{\M}{\mathcal{M}}

\newcommand{\Zd}{\mathbb{Z}^d}
\newcommand{\Rd}{\mathbb{R}^d}

\newcommand{\ep}{\varepsilon}

\newcommand{\K}{\mathcal{K}}

\renewcommand{\a}{\mathbf{a}}
\newcommand{\f}{\mathbf{f}}

\newcommand{\indc}{\mathds{1}}

\newcommand{\hslashslash}{%
  \raisebox{-.12ex}{%
    \scalebox{.6}{%
      \rotatebox[origin=c]{18}{$-$}%
    }%
  }%
}
\newcommand{\Lslash}{%
  {%
   \vphantom{L}%
   \ooalign{\kern-.08em\smash{\hslashslash}\hidewidth\cr$\,\!_{{L}}$\cr}%
  }%
}

\newcommand{\LL}{\Lslash}

\DeclareMathOperator{\dist}{dist}

\newcommand{\X}{\mathcal{X}}  
\newcommand{\Y}{\mathcal{Y}}

\renewcommand{\tilde}{\widetilde}

\begin{document}

\title{Calder\'on-Zygmund estimates for stochastic homogenization}

\begin{abstract}
We prove quenched~$L^p$--type estimates for the gradient of a solution of a quasilinear elliptic equation with random coefficients. 
\end{abstract}

\author[S. Armstrong]{Scott Armstrong}
\address{Ceremade (UMR CNRS 7534), Universit\'e Paris-Dauphine, Paris, France}
\email{armstrong@ceremade.dauphine.fr}

\author[J.-P. Daniel]{Jean-Paul Daniel}
\address{Institut f\"ur Geometrie, Fachrichtung Mathematik, Technische Universit\"at Dresden, 01062 Dresden, Germany}
\email{jean-paul.daniel@tu-dresden.de}

\date{\today}

\maketitle

\section{Introduction}

\subsection{Informal summary of results}

We study Calder\'on-Zygmund-type estimates for uniformly elliptic equations with random coefficients. We consider general quasilinear equations of the form
\begin{equation} \label{e.pde}
-\nabla \cdot  \left(\a \left( \nabla u,x \right)\right) = -\nabla \cdot \f \quad \mbox{in} \ B_R \subseteq \Rd.
\end{equation}
Here $\xi\mapsto \a(\xi,x)$ is a Lipschitz, uniformly monotone map, the maps $x\mapsto \a(\xi,x)$ are stationary random fields satisfying a finite range of dependence assumption, and the macroscopic scale~$R\gg1$ is typically large.

\smallskip

We are interested in obtaining~$L^p$-type bounds on~$|\nabla u|$, for large $p\in (2,\infty)$, which are independent of~$R$, in terms of the size of the vector field~$\f$ and the ellipticity of~$\a$. Such bounds generalize the classical Calder\'on-Zygmund estimate in the case~$\a(\xi,x) = \xi$ (when~\eqref{e.pde} is simply $-\Delta u = - \nabla \cdot \f$) which asserts that
\begin{equation} \label{e.CZest}
\left\| \nabla u \right\|_{L^p(B_{R/2})} \leq C(d,p) \left( \left\| \nabla u \right\|_{L^2(B_R)} + \left\| \f \right\|_{L^p(B_R)} \right).
\end{equation}
Recall that an estimate like~\eqref{e.CZest} fails to hold even in the linear case and even for smooth coefficients independently of~$R$, as a rescaling which maps $B_R$ to $B_1$ leads to coefficients which are rapidly oscillating. Indeed, the best available~$L^p$ estimate for gradients without assuming some control of the regularity of the coefficients is Meyers' estimate, which states that~\eqref{e.CZest} holds for every $p\in (2,2+\ep)$, for some tiny $\ep>0$ depending on the dimension and the ellipticity of the coefficients. 

\smallskip

Nevertheless, we show that an estimate similar to~\eqref{e.CZest} holds for equations with random coefficients satisfying mixing conditions. This can be considered a statistical effect, a manifestation of the principle that, at large scales, equations with random coefficients have better regularity properties than general equations because they are homogenizing. This idea originated in the work of Avellaneda and Lin~\cite{AL1} in the case of periodic coefficients and in~\cite{AL2} they proved Calder\'on-Zygmund-type estimates of the kind considered here.

\smallskip

Large-scale gradient estimates for solutions of elliptic equations with random coefficients was first proved in~\cite{AS}, which demonstrated that one can obtain regularity from homogenization in the stochastic setting, in the spirit of Avellaneda and Lin, by replacing the compactness methods of~\cite{AL1} with a quantitative approach and using a Campanato-type $C^{1,\alpha}$ iteration. Later, this approach was extended to handle more general equations and under relaxed mixing conditions in~\cite{GNO3} and \cite{AM}.

\smallskip

In particular,~\cite{AM} contains an $L^\infty$-type estimate for $|\nabla u|$ if $\f=0$ in~\eqref{e.pde}, or, more generally, if $\f$ is bounded and also random, satisfying a quantitative mixing condition. This estimate is the starting point for the results in this paper: we combine it with the ideas of Caffarelli and Peral~\cite{CP}, who introduced a general method for obtaining $L^p$ bounds from pointwise $L^\infty$-type bounds using the Calder\'on-Zygmund decomposition. The main source of difficulty here is that we do not have \emph{uniform} pointwise bounds, due to the randomness: the environment will have ``glitches" (regions where the pointwise estimate may fail, or local constant in the estimates is very large) and we must adapt the arguments to handle them. In Proposition~\ref{p.pert}, we formalize a flexible and modular tool which connects pointwise gradient estimates for~\eqref{e.pde} when $\f=0$ to gradient $L^p$ estimates. We believe it will have a wider applicability than the particular application here. 

\smallskip

Gradient estimates have played a central role in the theory of quantitative stochastic homogenization since the groundbreaking work of Gloria and Otto~\cite{GO1,GO2}, who demonstrated that they are the key to obtaining optimal scalings for the error, the fluctuations of the energy, and other  quantities of fundamental interest. One motivation behind the work in the present paper is to bound the error in the two-scale expansion in homogenization, which satisfies an equation with right-side in divergence form, such as~\eqref{e.pde}, but where the residual term $\f$ has a typical size much smaller than its $L^\infty$ norm (in other words, roughly speaking, $\| \f \|_{L^p} \ll \| \f \|_{L^\infty}$).  

\subsection{Statement of the main result}

Throughout the paper, we fix~$\Lambda \geq 1$ and an ambient dimension $d\geq2$.
We consider coefficient fields $$\a \in  L^\infty_{\mathrm{loc}}( \Rd \times\Rd ; \Rd)$$ satisfying, for every $\xi,\eta,x\in\Rd$,
\begin{equation} \label{e.a00}
\a(0,x) = 0,
\end{equation}
\begin{equation} \label{e.alip}
\left| \a(\xi,x) - \a(\eta,x) \right| \leq \Lambda \left| \xi-\eta \right|
\end{equation}
and 
\begin{equation} \label{e.aum}
\left( \xi-\eta \right) \cdot \left(\a(\xi,x) - \a(\eta,x) \right)\geq \frac1\Lambda \left| \xi-\eta \right|^2.
\end{equation}
We consider the set $\Omega$ of all such coefficient fields:
\begin{equation*} \label{}
\Omega:= \left\{ \a \in  L^\infty_{\mathrm{loc}}( \Rd \times\Rd ; \Rd) \,:\, \mbox{$\a$ satisfies~\eqref{e.a00},~\eqref{e.alip} and~\eqref{e.aum}}  \right\}.
\end{equation*}
We endow $\Omega$ with a family of $\sigma$--algebras $\{ \F_U \,:\, U \subseteq \Rd \ \mbox{is Borel} \}$ defined by 
\begin{equation*} \label{}
\F_U := \mbox{$\sigma$--algebra generated by } \quad \a \mapsto \int_{U} \a(\xi,x) \psi(x)\,dx,\quad \psi\in C^\infty_c(\Rd), \, \xi\in\Rd.
\end{equation*}
We drop the subscript if $U=\Rd$ and simply write $\F= \F_{\Rd}$. We denote the translation action of $\Rd$ on $\Omega$ by $\{ T_x \}_{x\in\Rd}$, that is, for each $y\in\Rd$, $T_y:\Omega\to \Omega$ denotes the map defined by
\begin{equation*}
(T_y\a)(\xi,x) : = \a(\xi,x+y).
\end{equation*}
We consider a probability measure~$\P$ on $(\Omega,\F)$ which is assumed to satisfy the following two conditions:

\begin{enumerate}

\smallskip

\item[(P1)] $\P$ is stationary with respect to $\Z^{d}$--translations: for every $z\in \Zd$ and $A\in \F$,
\begin{equation*} 
\P \left[ A \right] = \P \left[ T_z A \right].
\end{equation*}

\smallskip

\item[(P2)] $\P$ has a unit range of dependence: for all Borel subsets $U,V\subseteq \R^{d}$  such that $\dist(U,V) \geq 1$, we have that
\begin{equation*}
\F_U \quad \mbox{and} \quad \F_V \quad \mbox{are $\P$--independent.}
\end{equation*}
\end{enumerate}
Here we denote $\dist(U,V):= \inf\{ |x-y|\,:\, x\in U, \, y\in V \}$.

\smallskip

It is natural to state the main result in terms of ``coarsened" Lebesgue norms, because the regularizing effect we wish to observe happens only at large scales. These coarsened norms are designed to be blind to the behavior of functions at smaller scales. We define for $h>0$, $s\in [1,\infty)$, $\varphi \in L^1_{\mathrm{loc}}(\Rd)$ and $U\subseteq \Rd$:
\begin{align*} \label{}
\| \varphi \|_{L^{s}_h(U)} & := \left( \int_U \left( \fint_{B_h(x)} \left| \varphi(y) \right| \, dy\right)^s\, dx \right)^{\frac1s}, \\
\| \varphi \|_{\LL^{s}_{h}(U)} & := \left( \fint_U \left( \fint_{B_h(x)} \left| \varphi(y) \right| \, dy\right)^s\, dx \right)^{\frac1s} = |U|^{-\frac1s} \| \varphi \|_{L^{s}_h(U)}, \\
\| \varphi \|_{\LL^{s}(U)} & := \left( \fint_U \left| \varphi(x) \right|^s \, dx \right)^{\frac1s} =  |U|^{-\frac1s} \| \varphi \|_{L^{s}(U)}. 
\end{align*}
We introduce a ``coarsened" maximal operator $\M_h$ defined for each $h\geq 0$ and $\varphi \in L^1_{\mathrm{loc}}(\Rd)$ by
\begin{equation} \label{e.coarsemax}
\mathcal{M}_h(\varphi)(x) := \sup_{r>h}  \fint_{B_r(x)} |\varphi(y)| dy.
\end{equation}
If $\varphi \in L^1_{\mathrm{loc}}(U)$ for $U \subseteq \Rd$, then we define $\mathcal{M}_h(\varphi)(x)$ as above after extending the domain of $\varphi$ to $\Rd$ by taking $\varphi = 0$ in $\Rd \setminus U$. The usual maximal function is denoted by $\M=\M_0$.

\begin{theorem}[quenched $W^{1,p}$-type estimate]
\label{t.w1p}
Fix $2 < m < p$, $s\in \left(0,\frac{4}{m+2}\right)$ and $R\geq 10$. There exist $C(m,p,s,d,\Lambda) \geq 1$, $k(m,p,s,d)\geq 1$ and a nonnegative random variable $\Y_R$, depending on $(R,m,p,s,d,\Lambda)$, which satisfies 
\begin{equation} \label{e.intX}
\E \left[ \exp\left( \Y^s_R \right) \right] \leq C
\end{equation}
and such that the following holds: for every~$\a\in\Omega$, $\f \in L^2(B_R\,;\Rd)$, and solution $u\in H^1(B_R)$ in $B_R$ of the equation 
\begin{equation*} \label{}
-\nabla \cdot \left( \a\left( \nabla u,x \right) \right) = \nabla \cdot \f \quad \mbox{in} \ B_R,
\end{equation*}
we have the estimate
\begin{equation}
\label{e.w1p}
\left\| |\nabla u|^2 \right\|_{\LL^{m/2}_{1}(B_{R/2})} \\
 \leq \Y_R^2 M^2\left( \log (2+M) \right)^{k},
\end{equation}
where we define 
\begin{equation*}
M:=\left( \left\| \left| \nabla u \right|^2 \right\|_{L^1(B_{R})} + \left\| \left| \f\right|^2 \right\|_{\LL^{p/2}_{1}(B_{R})}   \right)^{\frac12}.
\end{equation*}
\end{theorem}

The fact that we must take $m< p$ in Theorem~\ref{t.w1p} is not an artifact of our method: it is an genuine effect of randomness. Indeed, since the random variable~$\X$ in Theorem~\ref{t.lip} is not bounded almost surely and considering a vector field~$\f$ whose support is concentrated on the set where $T_x\X$ is large, we expect any quenched $L^p$ estimate to give up some of the exponent.

\smallskip

The reason that the right side of~\eqref{e.w1p} has a logarithmic correction to its quadratic dependence in $M$ has to do with the interaction between the randomness of the coefficients and the fact the equation is nonlinear in general (see~\cite[Remark 1.2]{AM} and the remarks after the statement of Theorem~\ref{t.lip} in the next section). It cannot be scaled away because the assumption~(P2) has already fixed a length scale. It can however be removed if $\a$ is positively homogeneous in its first variable (e.g., if the equation is linear), and the right side of~\eqref{e.w1p} can be reduced to $\Y_R^2 M^2$.

\section{Proof of Theorem~\ref{t.w1p}}

As mentioned in the introduction, the main ingredient in the proof of the main result coming from the theory of homogenization is a quenched, gradient $L^\infty$-type estimate. For this purpose we use~\cite[Theorem 1.1]{AM}, as it is the most general and possesses optimal stochastic integrability, although a similar estimate under various other sets of assumptions would also do.

\begin{theorem}[{\cite[Theorem 1.1]{AM}}]
\label{t.lip}
Assume that $\P$ is a probability measure on $(\Omega,\F)$ and satisfies~(P1) and~(P2). Then for every $s\in (0,d)$, there exist a random variable~$\X\ge1$ and a constant~$C$, both of which depend on $(d,\Lambda,s)$, such that 
\begin{equation}
\label{e.lipintegrability}
\E \left[ \exp\left( \X^s \right) \right] \leq C
\end{equation}
and the following statement holds:  for every $R\geq10$ and $u\in H^1(B_R)$ satisfying
\begin{equation*}
-\nabla \cdot  \a(\nabla u,x)  = 0 \quad \mbox{in} \ B_R,
\end{equation*}
and setting
\begin{equation*}
M:= \frac1R \inf_{a\in\R} \left( \fint_{B_R} \left| u(x) - a \right| \,dx \right)^{\frac12},
\end{equation*}
we have the estimate
\begin{equation}
\label{e.Linfty}
\fint_{B_r} \left| \nabla u(x) \right|^2 \,dx \leq CM^2 \quad \mbox{for every} \quad \X \log (2+M) \leq r \leq \frac12 R.
\end{equation}
\end{theorem}

As explained in~\cite{AM}, we consider~\eqref{e.Linfty} to be an $L^\infty$--type bound because the radius~$r$ can be taken to be on the order of the microscopic scale (here of order one), while the macroscopic scale, given by $R\gg1$, is much larger. Indeed, notice that ~\eqref{e.Linfty} can be written in terms of the coarsened maximal operator (defined in~\eqref{e.coarsemax}) as 
\begin{equation*}
\M_{r_*} \left( \left| \nabla u \right|^2 \right)(0) \leq CM^2,
\end{equation*}
where $r_*:= \X \log (1+M)$ denotes the ``minimal radius." In particular, we deduce 
\begin{equation}
\label{e.pointwiseform}
\M_1\left( \left| \nabla u \right|^2 \right)(0) \leq Cr_*^{d}M^2 = C \X^d M^2 \log^d (2+M),
\end{equation}
so $\X$ can also be thought of as a (random) constant in an {\it a priori} local gradient bound which is blind to oscillations on scales smaller than one. 

\smallskip

The reason that the minimal radius $r_*$  has some weak dependence on $M$ has to do with the interaction between the randomness of the coefficients and the fact the equation is nonlinear (see~\cite[Remark 1.2]{AM}). It can be removed completely in the case that $\a$ is positively homogeneous in its first variable (thus in particular when the equation is linear) and we have simply $r_* := \X$. 

\smallskip

We next turn to the link between Theorems~\ref{t.lip} and~\ref{t.w1p}, which is formalized in the following proposition. The proof appears in the next section and is a modified Calder\'on-Zygmund-type argument in the spirit of Caffarelli and Peral~\cite{CP}. 

\smallskip

Motivated by Theorem~\ref{t.lip}, especially with its conclusion written in the form~\eqref{e.pointwiseform}, we define, for a given nonnegative measurable function $\K:\Rd \to [0,\infty)$ and a ball $B_{2r}(x_0)$ with radius $2r\geq 2$, the set 
\begin{multline*} \label{}
\mathcal A(B_{2r}(x_0),\K) := \bigg\{ v \in H^1(B_{2r}(x_0)) \, :  \,  \mbox{for every} \  x\in B_r(x_0), \\
\fint_{B_1(x)} \left| \nabla v(y)\right|^2\,dy  
\leq \K(x) \fint_{B_{2r}(x_0)} \left| \nabla v(y) \right|^2\,dy \bigg\}.
\end{multline*}
Roughly speaking, this is the collection of functions satisfying a gradient $L^\infty$-type estimate, similar to~\eqref{e.Linfty}, where the quality of the estimate varies and is given by~$\mathcal K$.

\smallskip

\begin{proposition}
\label{p.pert}
Fix $2 < p < q < \infty$, $R\geq 10$, $f\in L^2(B_R)$, $u\in H^1(B_R)$ and a nonnegative measurable function $\K:\Rd \to[1,\infty)$. Assume that~$u$ has the following property: for every $r\in [4,R/2]$ and $x_0 \in B_{R/2}$, 
\begin{equation*} \label{}
\inf_{v\in \mathcal A\left(B_{\frac r2}(x_0),\,\K\right)} \, \fint_{B_{\frac r2}(x_0)} \left| \nabla u(x) - \nabla v(x) \right|^2 \, dx 
\leq  \fint_{B_r(x_0)} \left|f(x) \right|^2 \, dx.
\end{equation*}
Then there exists~$C(d,p,q) \geq 1$ such that, for every $t>0$,
\begin{equation*}
\left| \left\{ x \in B_{R/2} : \mathcal M_1(|\nabla u|^2) (x) > t \right\} \right| 
\leq C \left|B_{R/2}\right| \left(\frac{t}{t_\ast}\right)^{- \frac{p}{2}\left(1 - \theta \right)}  \left \| \K^{q/2} \right\|_{\LL^{1}_1(B_R)}^{\theta},
\end{equation*}
where the exponent $\theta\in (0,1)$ is given by
\begin{equation} \label{e.theta}
\theta := \frac{p^2+2p}{p^2+2q}
\end{equation}
and $t_\ast \in \R$ is given by
\begin{equation*} 
t_\ast:= \left \| |\nabla u|^{2} \right\|_{\LL^{1}_1(B_R)}+\left \|f^{2} \right\|_{\LL^{p/2}_1(B_R)}.
\end{equation*}
\end{proposition}

We next present the proof of the main result of the paper, which is obtained by combining Theorem~\ref{t.lip} and Proposition~\ref{p.pert}. The proof of the proposition is given in the following section.

\begin{proof}[{Proof of Theorem~\ref{t.w1p}}]
Fix~$2< m< p$, $s\in (0,4/(m+2))$ and $R\geq 10$. 
Throughout this argument, we let~$C$ and~$c$ denote positive constants depending only on $(d,\Lambda,m,p)$ and which may vary in each occurrence. By taking $p$ to be closer to $m$, if necessary, by an amount depending only on $s$, we may suppose that $s < 4/(p+2)-c$. We also select an exponent $q>p$, depending only on $(m,p,s)$, such that $q\leq C$, 
\begin{equation}
\label{e.scont1}
s  \leq \frac{4(q-p)}{(p+2)q} - c
\end{equation}
and
\begin{equation*}
m \leq p(1-\theta)- c,
\end{equation*}
where, as in the statement of Proposition~\ref{p.pert}, we denote $\theta:=(p^2+2p)/(p^2+2q)$. 
Finally, we set
\begin{equation} \label{e.ndef}
n:=\max\{ \tfrac{1}{2}, s(p+2)q/4(q-p) \}
\end{equation}
and observe by~\eqref{e.scont1} that $n\leq 1-c$. Let~$\X$ denote the random variable in the statement of Theorem~\ref{t.lip} such that~\eqref{e.lipintegrability} holds with exponent~$s\in (dn+c,d)$. 

\smallskip

We let $\K$ denote the random variable $\X^{d}\log^d (2+\X)$. We may extend $\X$ and $\K$ to be a $\Zd$--stationary random field on~$\Rd$ by setting $\X(x): = T_z\Z$ and $\K(x):= T_x\K$ for every $x\in\Rd$. It follows from~\eqref{e.lipintegrability} that, for every $x\in \Rd$, 
\begin{equation} \label{e.Kintegrx}
\sup_{x\in \Rd} \E \left[ \exp\left( \left| \K(x) \right|^n\right) \right] \leq C.
\end{equation}
Indeed, this holds for every $x\in\Zd$ immediately from~\eqref{e.lipintegrability} by~$\Zd$--stationarity. We see that it also holds for every~$x\in\Rd$ by applying Theorem~\ref{t.lip} to the pushforward of~$\P$ under the map~$L\mapsto T_x\a$, which satisfies the same assumptions as~$\P$.

\smallskip

\emph{Step 1.} We show that, for every $r\in [4,R/2]$ and $x_0\in B_{R/2}$, 
\begin{equation}
\label{e.verapprox}
\inf_{v\in \mathcal A\left( B_{\frac r2}(x_0),C\K \log^{d}(2+M)\right) } \fint_{B_{\frac r2}(x_0)} \left| \nabla u(x) - \nabla v(x) \right|^2 \, dx 
\leq C\fint_{B_r(x_0)}  \left| \f(x)\right|^2 \,dx .
\end{equation}
Let $v \in u + H^1_0(B_r(x_0))$ be the unique solution of the Dirichlet problem 
\begin{equation*}
\left\{ 
\begin{aligned}
& -\nabla \cdot \left( \a(\nabla v,x) \right) = 0 & \mbox{in} & \ B_r(x_0), \\
& v = u &  \mbox{on} &\  \partial B_r(x_0).
\end{aligned}
\right.
\end{equation*}
By Caccioppoli's inequality, 
\begin{equation*}
\fint_{B_{\frac r2}(x_0)} \left| \nabla u(x) - \nabla v(x) \right|^2 \, dx 
\leq C\fint_{B_r(x_0)}  \left| \f(x)\right|^2 \,dx.
\end{equation*}
By Theorem~\ref{t.lip}, we have that $v\in \mathcal A\left( B_{\frac r2}(x_0), \tilde K \right)$, where 
\begin{align*}
\tilde K(x) & = C \X^d(x) \log^d \left( 2+ \fint_{B_r(x_0)} \left| \nabla u(x) \right|^2\,dx \right)  \\
& \leq C\X^{d} (x)  \log^{d}(2+M)  \log^d \left( 2+ \X(x) \right) = C\K(x) \log^{d}(2+M). 
\end{align*}
%
Here we used the fact that, for any $x$ and any ball $B_r(x_0)$ containing $x$ with $r\geq 1$, we have
\begin{equation*}
 \fint_{B_r(x_0)} \left| \varphi(x) \right|\,dx \leq C \M_1(\varphi)(x).
\end{equation*}
This completes the proof of~\eqref{e.verapprox}.

\smallskip

\emph{Step 2.} We apply Proposition~\ref{p.pert}. We obtain, for $t_\ast:= M^2$ and every $t>C t_\ast$, 
\begin{multline*} 
|B_{R/2}|^{-1}  \left| \left\{ x \in B_{R/2} : \mathcal M_1(|\nabla u|^2) (x) > t \right\} \right| 
\\
\leq C\log^{d\theta q/2}(2+M) \left(\frac{t}{t_\ast} \right)^{- \frac{p}{2}\left(1 - \theta \right)}
\left \| \K^{q/2} \right\|_{\LL^{1}_1(B_R)}^{\theta}.
\end{multline*}
By integrating the previous inequality, using that $m \leq p(1-\theta) - c$, we obtain, 
for every $\tau \geq Ct_\ast$, 
\begin{align*}
\lefteqn{ 
\fint_{B_{R/2}} \left| \M_1\!\left(|\nabla u|^2 \right)\!(x) \right|^{\frac m2} \,dx  
} \qquad & \\
& = |B_{R/2}|^{-1} \int_{0}^\infty t^{\frac m2-1} \left| \left\{ x \in B_{R/2} : \mathcal M_1(|\nabla u|^2) (x) > t \right\} \right|\,dt \\
& \leq \int_0^\tau t^{\frac m2-1} \,dt + t_\ast^{\frac p2 (1 - \theta)} \log^{d\theta q/2}(2+M) 
 \left \| \K^{q/2} \right\|_{\LL^{1}_1(B_R)}^{\theta}  \int_\tau^\infty t^{\frac m2-1-\frac p2(1-\theta)} \, dt \\
 & \leq C\tau^{\frac m2} \left( 1+ 
 \left(\frac{\tau}{t_\ast} \right)^{-\frac p2 (1 - \theta)}
\log^{d\theta q/2}(2+M) 
 \left \| \K^{q/2} \right\|_{\LL^{1}_1(B_R)}^{\theta} 
\right).
\end{align*}
Noticing that $2\theta/(p(1-\theta)) = (p+2)/(q-p)$ and substituting
\begin{equation*}
\tau: = 
C t_\ast  \left(  \log^{dq/2}(2+M)  \left \| \K^{q/2} \right\|_{\LL^{1}_{1}(B_R)} \right)^{\frac {p+2}{q- p}}  
\end{equation*}
into the previous inequality yields
\begin{equation*}
\left\|  \M_1\!\left(|\nabla u|^2 \right) \right\|_{\LL^{m/2}(B_{R/2})} 
\leq C \Y_R^2  M^2 \left( \log (2+M) \right)^{\frac{dq(p+2)}{q-p}},
\end{equation*}
where we  defined $\Y_R$ to be the random variable  
\begin{equation*}
\Y_R:= C \left(\fint_{B_R}  \left| \K(x) \right|^{\frac{q}{2}} \, dx\right)^{\frac{p+2}{2(q- p)}}.
\end{equation*}
The previous inequality is a stronger form of the desired inequality~\eqref{e.w1p}, since
\begin{equation*} \label{}
\left\| |\nabla u|^2 \right\|_{\LL^{m/2}_{1}(B_{R/2})} \leq \| \M_1(|\nabla u|^2)(x)\|_{\LL^{m/2}(B_{R/2})}.
\end{equation*}

\smallskip

\emph{Step 3.} We complete the proof by verifying that~$\Y_R$ satisfies~\eqref{e.intX}. Define
\begin{equation*}
\mathcal{Z}_R:= \fint_{B_R} \left| \K(x) \right|^{\frac q2}\,dx.
\end{equation*}
With $n$ defined in~\eqref{e.ndef}, we claim that
\begin{equation*} 
\E \left[ \exp\left( \mathcal{Z}_R^{{2n}/{q}} \right)\right]  \leq C. 
\end{equation*}
This is a consequence of~\eqref{e.Kintegrx} and Jensen's inequality:
\begin{align*}
\E \left[ \exp\left( \mathcal{Z}_R^{{2n}/{q}} \right)\right] & 
=  \E \left[ \exp\left( \left( \fint_{B_R} \left| \K(x) \right|^{\frac q2}\,dx \right)^{{2n}/{q}} \right)\right] \\
& \leq \E \left[ \fint_{B_R} \exp\left( \left| \K(x) \right|^n \right) \,dx  \right] + C\\
& = \fint_{B_R} \E \left[  \exp\left( \left| \K(x) \right|^n \right)   \right] \,dx +C \\
& \leq C.
\end{align*}
To get the second line in the previous inequality string, we used Jensen's inequality and the fact that (since $2n/q \geq q \geq c$) the map~$t\mapsto \exp\left(t^{2n/q}\right)$ is convex on the interval~$[C,\infty)$. Here is a more detailed derivation:
\begin{align*}
\exp\left( \left( \fint_{B_R} \left| \K(x) \right|^{\frac q2}\,dx \right)^{{2n}/{q}} \right) & \leq  \exp\left( \left( \fint_{B_R} \max\left\{ C, \left| \K(x) \right|^{\frac q2}\right\} \,dx \right)^{{2n}/{q}} \right) \\
& \leq \fint_{B_R} \exp \left(   \max\left\{ C, \left| \K(x) \right|^{\frac q2}\right\}^{2n/q}  \right) \, dx \\
& \leq \fint_{B_R} \left( \exp\left(\left| \K(x) \right|^{n} \right) +\exp\left( C^{2n/q} \right)\right)  \,dx \\
& \leq \fint_{B_R}\exp\left(\left| \K(x) \right|^{n} \right)\,dx  + C.
\end{align*}
Using that $s \leq 4n(q-p)/(p+2)q$ and rewriting the above inequality for $\mathcal{Z}_R$ in terms of~$\Y_R$, we get
\begin{equation*}
\E \left[ \exp\left( \Y_R^{s} \right)\right] \leq \E \left[ \exp\left( \Y_R^{{4n(q-p)}/{(p+2)q}} \right)\right]  =  \E \left[ \exp\left(\mathcal{Z}_R^{2n/q} \right) \right] \leq C.
\end{equation*}
This completes the proof of the theorem. 
\end{proof}

\section{The Calder\'on-Zygmund argument}
\label{ss.pert}

As mentioned above, the proof of Proposition~\ref{p.pert} is a modification of the Calder\'on-Zygmund method introduced in~\cite{CP}. We structure the argument somewhat differently, however: similar to~Byun~\cite{B} and in contrast to~\cite{CP}, we use the Vitali covering theorem rather than a cube decomposition for our covering theorem needs. We remark that it is possible to avoid the use of maximal functions by using an argument like the one in Acerbi and Mingione~\cite{AMing}. A statement like Proposition~\ref{p.pert} in the case that $\K$ is bounded has been previously proved by Shen~\cite{Shen}.

\smallskip

We begin with some elementary observations concerning the coarsened maximal function defined in~\eqref{e.coarsemax}, above. The proofs of these facts are omitted, since they are easy modifications of the arguments for $\M_0$ which are classical. First, notice by the definition that $\M_s \leq \M_t$ if $t\leq s$. For every~$\varphi \in L^1(\Rd)$, we have the weak-type~$L^1$ estimate
\begin{equation} \label{e.MweakL1}
\left| \left\{ x \in \Rd : \M_h (\varphi) (x)> t \right\} \right| 
\leq \frac{C(d)}{t}  \left\|  \varphi \right\|_{L^{1}_h(\Rd)} , 
\end{equation}
 and, for $p>1$ and $\varphi\in L^1(\Rd)$, the strong-type~$L^p$ estimate
\begin{equation} \label{e.MstrongLp}
  \left\|  \M_h(\varphi) \right\|_{L^p(\Rd)}^p
\leq C(p,d)  \left\|  |\varphi|^2 \right\|_{L^{p/2}_h(\Rd)}^{p/2} .
 \end{equation}
We also recall the classical statement that 
\begin{equation} \label{e.infMball}
\inf_{x\in B_{R/8}} \M (\varphi)(x) \leq C(d) \fint_{B_R} |\varphi(x)|\,dx.
\end{equation}

\smallskip

\begin{proof}[{Proof of Proposition~\ref{p.pert}}]
Throughout the argument, we let $C$ and $c$ denote positive constants which depend only on $(p,q,d)$ and may vary in each occurrence. We denote the sublevel sets of $\M_1(|\nabla u|^2)$ by
\begin{equation*}
A(t) : = \left\{ x \in B_R : \mathcal M_1(|\nabla u|^2) (x) \leq t \right\}.
\end{equation*}
The main step in the argument is to prove that, for every $\sigma \in (0,1]$, $\omega,t > 0$, $x_0\in B_{R/2}$ and $r\in [1,R/8]$ satisfying 
\begin{equation}\label{e.omegsig}
\omega^{q/(q-2)}\sigma^{-2/(q-2)}  \geq C,
\end{equation}
\begin{equation} \label{e.contact}
 B_r(x_0) \cap A(t)\neq \emptyset,
\end{equation}
\begin{equation} \label{e.intcondf}
 \fint_{B_{4r}(x_0)} \left| f(x) \right|^2\,dx     \leq  \sigma t 
\end{equation}
and
\begin{equation} \label{e.intcondK}
 \fint_{B_r(x_0)} \left|\K(x)\right|^{\frac q2}\,dx \leq \omega^{\frac q2},
\end{equation}
we have 
\begin{equation} \label{e.climbs}
\left| B_r(x_0)  \setminus A\left( \omega^{\frac q{q-2}} \sigma^{-\frac 2{q-2}}  t\right) \right| \leq C  \sigma |B_r| \cdot \omega^{-\frac q{q-2}} \sigma^{\frac 2{q-2}}.
\end{equation}
The proof of this statement is accomplished in the next step, and then in Steps~2 and~3 we use it to complete the proof of the proposition. 

\smallskip

\emph{Step 1.} We fix $\sigma\in (0,1]$, $\omega,t >0$, $x_0 \in B_{R/2}$ and $r \in [1,R/8]$ such that~\eqref{e.omegsig}, \eqref{e.contact}, \eqref{e.intcondf} and~\eqref{e.intcondK} hold and proceed to derive~\eqref{e.climbs}.  

\smallskip

Observe that~\eqref{e.contact} implies that
\begin{equation} \label{e.sillycovbnd}
\sup_{s\geq r} \fint_{B_{s}(x_0)} \left|\nabla u(x) \right|^2\, dx \leq 2^dt.
\end{equation}
By the hypotheses of the proposition, there exists 
\begin{equation} \label{e.vinAinfty}
v\in  \mathcal A(B_{2r}(x_0),\K) 
\end{equation}
such that
\begin{equation} \label{e.approxubyv}
\fint_{B_{2r}(x_0)} \left| \nabla u(x) - \nabla v(x) \right|^2 \, dx
  \leq \fint_{B_{4r}(x_0)} \left|f(x) \right|^2 \, dx .
\end{equation}
Note that~\eqref{e.intcondf},~\eqref{e.sillycovbnd},~\eqref{e.approxubyv} and $\sigma \leq 1$ imply that 
\begin{equation} \label{e.supsillyvbnv}
\fint_{B_{2r}(x_0)} \left|\nabla v(x) \right|^2\, dx \leq Ct + \fint_{B_{4r}(x_0)} \left|f(x) \right|^2 \, dx
   \leq Ct. 
\end{equation}
By~\eqref{e.sillycovbnd}, we have, for every $x \in B_r(x_0)$,
\begin{align*} \label{}
\M_1\!\left(\left| \nabla u\right|^2 \right)\!(x) 
& \leq 2 \M_1\!\left(\left| \nabla v\right|^2 \indc_{B_{r}(x_0)}\right)(x) + 2\M_1\!\left(\left| \nabla u - \nabla v\right|^2\indc_{B_{r}(x_0)} \right)(x)+ Ct.
\end{align*}
We deduce that, for every $s\geq C$,
\begin{multline*}
\left| \left\{ x \in B_{r}(x_0) \,:\, \M_1\!\left(\left| \nabla u\right|^2 \right)(x) > st \right\}  \right| \\
 \leq \left| \left\{ x \in B_{r}(x_0) \,:\, \M_1\!\left(\left| \nabla v\right|^2 \indc_{B_{r}(x_0)}\right)(x) > \tfrac14st \right\} \right| \\
+ \left| \left\{ x \in B_{r}(x_0) \,:\, \M_1\!\left(\left| \nabla u-\nabla v\right|^2 \indc_{B_{r}(x_0)}\right)(x) > \tfrac14st \right\} \right|.
\end{multline*}
The first term on the right side is controlled by Chebyshev's inequality,~\eqref{e.MstrongLp}, \eqref{e.intcondK}, \eqref{e.vinAinfty}, \eqref{e.supsillyvbnv} and~$r\geq 1$ and the second term by~\eqref{e.MweakL1},~\eqref{e.intcondf} and~\eqref{e.approxubyv}, as follows:
\begin{align*}
\lefteqn{ \left| \left\{ x \in B_{r}(x_0) \,:\, \M_1\!\left(\left| \nabla v\right|^2\indc_{B_{r}(x_0)} \right)(x) > \tfrac14st \right\} \right|} \qquad & \\
&  \leq C(st)^{-\frac q2} \left(\int_{ B_{r}(x_0)} \left| \M_1\!\left(\left| \nabla v\right|^2\indc_{B_{r}(x_0)} \right)(x) \right|^{\frac q2}\,dx \right)\\
& \leq C(st)^{-\frac q2} \int_{\Rd} \left( \fint_{B_1(x)} \left|\nabla v(y) \right|^2 \indc_{B_r(x_0)}(y) \,dy\right)^{\frac q2}\,dx  \\
& \leq C|B_r| (st)^{-\frac q2} \left( \fint_{B_r(x_0)} \left|\K(x)\right|^{\frac q2}\,dx \right)  
\left( \fint_{B_{2r}(x_0)} \left| \nabla v(x) \right|^2\,dx \right)^{\frac q2} \\
& \leq C|B_r|s^{-\frac q2} \omega^{\frac q2}.
\end{align*}
and
\begin{multline*}
\left| \left\{ x \in B_{r}(x_0) \,:\, \M_1\!\left(\left| \nabla u-\nabla v\right|^2 \indc_{B_{r}(x_0)} \right)(x) > \tfrac14st \right\} \right| \\
 \leq \frac{C}{st}\int_{B_{2r}(x_0)} \left| \nabla u(x)-\nabla v(x)\right|^2 \, dx \leq \frac{C}{st}\int_{B_{4r}(x_0)} \left|f(x) \right|^2 \, dx \leq \frac{C\sigma }{s}  |B_r|. 
\end{multline*}
Taking $s:= \omega^{q/(q-2)}\sigma^{-2/(q-2)}$, using~\eqref{e.omegsig} and combining the above yields, 
\begin{equation*} \label{}
\left| \left\{ x \in B_{r}(x_0) \,:\, \M_1\!\left(\left| \nabla u\right|^2 \right)(x) > \left(\omega^{\frac q{q-2}} \sigma^{-\frac 2{q-2}} \right) t \right\}  \right| \leq  C\sigma |B_r| \cdot \omega^{-\frac q{q-2}} \sigma^{\frac 2{q-2}}.
\end{equation*}
This is~\eqref{e.climbs}.

\smallskip

\emph{Step 2.} According to~\eqref{e.infMball}, if $C_0 \geq C$, then the parameter
\begin{equation*} 
 t_0 := C_0 t_\ast =  C_0 \left( \left \| |\nabla u|^2 \right\|_{\LL^1_1(B_R)} + \left \| f^2 \right\|_{\LL_1^{p/2}(B_R)}
 \right)
\end{equation*}
satisfies, for every $x\in B_{R/2}$, 
\begin{equation*} 
B_{R/40}(x) \cap A(t_0) \neq \emptyset. 
\end{equation*}
We show next that, for every $\sigma \in (0,1]$, $\omega>0$ satisfying~\eqref{e.omegsig} and $t\geq t_0$, 
\begin{align} \label{e.vitaliclm}
\lefteqn{ \bigg| B_{R/2}  \setminus   A\big( \omega^{\frac q{q-2}}   \sigma^{-\frac 2{q-2}}  t\big) \bigg| } \qquad & \\ 
&  \begin{multlined}[.8\textwidth]
\leq C \omega^{-\frac q{q-2}} \sigma^{\frac q{q-2}}  \left| B_{R/2} \setminus A(t)\right| + C \left| \left\{ x\in B_{R} : \M_1\!\left( |f|^2\right) (x) > c \sigma t \right\} \right|  \\
  +  C\big| \big\{ x\in B_{R} : \M_1( |\K|^{\frac q2}) (x) > c\omega^{\frac q2} 
 \big\} \big|.
\end{multlined}\nonumber
\end{align}
We first notice that, for every $x\in B_{R/2}$,
\begin{equation} \label{e.bndryeats}
\dist(x,A(t)) \leq 2 \quad \implies \quad x\in A(Ct).
\end{equation}
In preparation for the application of the Vitali covering theorem, we observe that, for every $x\in B_{R/2}$ such that $\dist(x,A(t))> 2$, there exist $r\in [1,R/40]$ and $y\in B_{R/2-1}$ such that $x\in B_r(y)$, $B_{r/2}(y) \subseteq B_{R/2} \setminus A(t)$ and $B_r(y) \cap A(t) \neq \emptyset$. 
Indeed, consider the ball $B_r(y_0)$ with $r$ slightly bigger than $r_0$ given by 
\begin{equation*}
r_0:=\min \big\{ r\geq 1 : B_{r}(y) \cap A(t) \neq \emptyset  
\text{ with }y= x -r x/|x|
\big\}  .
\end{equation*}
and $y_0= x -r_0 x/|x| $. In particular, observe  $r_0$ is well-defined and less than $R/40$, $B_{r/2}(y_0) \subseteq B_{3r_0/4}(y_0) \subseteq B_{R/2}\setminus A(t)$ and $x\in B_r(y)$.

\smallskip

The Vitali covering theorem yields a finite collection of pairwise disjoint balls $\{ B_{r_i} (y_i)\}_{i=1}^N$ such that $r_i \in [1,R/40]$, $B_{r_i/2}(y_i)\subseteq B_{R/2} \setminus A(t)$, $B_{r_i}(y_i) \cap A(t) \neq \emptyset$ and
\begin{equation} \label{e.squishy}
\left\{ x\in B_{R/2} \,:\, \dist(x,A(t))>2 \right\} \subseteq \bigcup_{i=1}^N B_{5r_i}(y_i).
\end{equation}
Put $\beta:= \omega^{q/(q-2)}\sigma^{-2/(q-2)}$. Applying the implication proved in Step~1 to each ball $B_{5r_i}(y_i)$, we deduce that, for each $i\in \{1,\ldots,N\}$, at least one of the following three alternatives must hold:
\begin{enumerate}
\item[(1)] $\displaystyle\left|B_{5r_i}(y_i) \setminus A(\beta t) \right|\leq \sigma\beta^{-1} \left|B_{5r_i}\right|$, \smallskip
\item[(2)]  $\displaystyle\fint_{B_{20r_i}(y_i)} |f(x)|^2\,dx 
> \sigma t$, 
\item[(3)] $\displaystyle\fint_{B_{5r_i}(y_i)} \left| \K(x) \right|^{\frac q2} \, dx > \omega^{\frac q2}.$
\end{enumerate}
We may partition $\{ 1,\ldots, N\}$ into disjoint subsets $I_1$, $I_2$ and $I_3$ such that the $j$th alternative above holds for each $i\in I_j$. Note in particular that
\begin{equation*}
\left\{ \begin{aligned} 
& i\in I_2 \quad \implies\quad \M_1(|f|^2)>c\sigma t \quad \mbox{in} \ B_{r_i/2}(y_i), \\
& i\in I_3 \quad \implies\quad \M_1(|\K|^{\frac q2})>c \omega^{\frac q2} \quad \mbox{in} \ B_{r_i/2}(y_i).
\end{aligned} \right.
\end{equation*}
Using~\eqref{e.bndryeats},~\eqref{e.squishy}  and $\beta \geq C$ (and enlarging $C$ if necessary) and that the balls $\{ B_{r_i/2}(y_i)\}$ are pairwise disjoint and contained in $B_{R/2} \setminus A(t)$, we thus obtain
\begin{align*} \label{}
\left|B_{R/2} \setminus A(\beta t) \right| & \leq \sum_{i=1}^N \left|B_{5r_i}(y_i) \setminus A(\beta t) \right| \\
& \leq C\sigma\beta^{-1} \sum_{i \in I_1} |B_{r_i/2}| + \sum_{i\in I_2 \cup I_3} |B_{r_i/2}| \\
& \begin{multlined}
 \leq C\sigma \beta^{-1} \left|B_{R/2} \setminus A(t) \right| + C \left| \left\{ x \in B_{R} \,:\, \M_1\!\left(|f|^2\right)(x) > c\sigma t \right\}\right|    \\
+ C \big| \big\{ x \in B_{R} \,:\,  \M_1(|\K|^{\frac q2})>c \omega^{\frac q2}    \big\}\big|.
\end{multlined}
\end{align*}
This is~\eqref{e.vitaliclm}.

\smallskip

\emph{Step 3.}
We complete the argument by iterating~\eqref{e.vitaliclm}. Fix $T \geq t_0=C_0t_\ast$.  Define
\begin{equation*} \label{}
\left\{ \begin{aligned}
\sigma & := \left( \frac{T}{t_\ast} \right)^{-\frac{p}{2m}}
\left \| \K^{q/2} \right\|_{\LL^{1}_1(B_R)}^{-\frac1m}, \\
\beta &: = c_0^{\frac2{p-2}}  \sigma^{- \frac2{p-2} }, \\
\omega &:= \beta^{\frac {q-2}q} \sigma^{ \frac 2q} ,
\end{aligned}
\right.
\end{equation*}
where $0< c_0\leq \frac12$ is chosen below, $\theta$ is defined in~\eqref{e.theta},  
\begin{equation*}
m  := \left(\frac{p}{2} + \frac{p+q}{p-2}\right) \quad \mbox{and} \quad  \nu : = \frac{p^2}{p^2+2q}.
\end{equation*}
Select $k\in\N$ to be such that $\beta^{k} t_0 < T \leq \beta^{k+1} t_0$. 
The choice of $t_0$ implies that
\begin{equation*}
\beta=  c_0^{\frac2{p-2}}\left(\frac{T}{t_\ast} \right)^{\frac{2\nu}{p}}  \left \| \K^{q/2} \right\|_{\LL^{1}_1(B_R)}^{\frac{4\nu}{p^2}}
\geq c_0^{\frac 2{p-2}} C_0^{\frac{2\nu}{p}},
\end{equation*}
where we used that $\K\geq 1$ and $T\geq t_0$.
Therefore the inequality~\eqref{e.vitaliclm} is valid provided that
\begin{equation} \label{e.c0C0}
c_0^{\frac 2{p-2}} C_0^{\frac{2\nu}{p}} \geq C.
\end{equation}
In this case we may iterate it to obtain 
\begin{align*}
\big| B_{R/2} \setminus  A& \left( \beta^kt_0\right) \big|  \leq C^k\left( \frac{\sigma}{\beta}\right)^k |B_{R/2}| \\
 &  + C\sum_{j=0}^{k-1} C^{k-1-j}\left( \frac{\sigma}{\beta}\right)^{k-1-j} \left|\left\{ x\in B_{R} : \M_1\!\left( |f|^2\right) (x)  > c \sigma \beta^j t_0 \right\} \right| \\
 & + C \sum_{j=0}^{k-1} C^{k-1-j} \left(\frac{\sigma}{\beta}\right)^{k-1-j} \big|\big\{ x\in B_{R} : \M_1( |\K|^{\frac q2}) (x) > c\omega^{\frac q2}  \big\} \big|.
\end{align*}
By making $C_0$ larger, if necessary, we may assume that $\sigma \leq c$, which implies that $\sigma \beta^{-1} \leq \sigma \leq c$, and thus we may simplify the last term in the previous inequality to obtain 
\begin{multline} \label{e.crabtcha}
 \left| B_{R/2}  \setminus A\left( \beta^kt_0 \right) \right|  \leq C^k\left( \frac{\sigma}{\beta}\right)^k |B_{R/2}|  
 + C \big|\big\{ x\in B_{R} : \M_1( |\K|^{\frac q2}) (x) > c\omega^{\frac q2}   \big\} \big|\\
  \qquad + C\sum_{j=0}^{k-1}C^{k-1-j} \left( \frac{\sigma}{\beta}\right)^{k-1-j} \left|\left\{ x\in B_{R} : \M_1\!\left( |f|^2\right) (x)  > c \sigma \beta^j t_0 \big\}  \right\} \right|.
 \end{multline}
We proceed by estimating the terms on the right side of~\eqref{e.crabtcha} each in turn. 

For the first term on the right  of~\eqref{e.crabtcha}, we use the fact that 
$\sigma \beta^{-1} = c_0 \beta^{-\frac {p}2} $ and the choice of~$k$ to obtain 
\begin{equation*} \label{}
C^k \left( \frac{\sigma}{\beta}\right)^k |B_{R/2}| = C^kc_0^{k}\left( \beta^{k} \right)^{-\frac p2}|B_{R/2}| \leq C^kc_0^{k} \beta^{\frac p2} \left(\frac{T}{t_0} \right)^{-\frac p2} |B_{R/2}|.  
\end{equation*}
We now fix $c_0 \in (0,\frac12]$ to be sufficiently small, depending on the appropriate quantities (and also fix $C_0$ larger, if necessary, so that~\eqref{e.c0C0} holds),
 to obtain
\begin{equation} \label{e.boink1}
C^k \left( \frac{\sigma}{\beta}\right)^k |B_{R/2}|  \leq |B_{R/2}| \beta^{\frac p2}\left(\frac{T}{t_0} \right)^{-\frac p2} \leq C |B_{R/2}| \left(\frac{T}{t_0} \right)^{-\frac p2+\nu}\left \| \K^{q/2} \right\|_{\LL^{1}_1(B_R)}^{\frac{2\nu}p}.
\end{equation}
The last inequality was obtained using the observation that
\begin{equation*} \label{}
\beta  = c\sigma^{-\frac 2{p-2}} = 
c\left(\frac{T}{t_0}\right)^{\frac{2\nu}p} 
\left \| \K^{q/2} \right\|_{\LL^{1}_1(B_R)}^{\frac{4\nu}{p^2}}.
\end{equation*}
For the second term on the right of~\eqref{e.crabtcha}, we observe that
\begin{align*} \label{}
\omega^{-\frac q2}  = c \sigma^{m \left( 1 - \theta \right)} 
 =c \left( \frac{T}{t_\ast}\right)^{-\frac p2\left( 1-\theta \right) }   \left \| \K^{q/2} \right\|_{\LL^{1}_1(B_R)}^{\theta-1}
\end{align*}
and thus
\begin{align}\label{e.boink2}
 \big|\big\{ x\in B_{R}:\M_1( |\K|^{\frac q2})(x) > c\omega^{\frac q2}  \big\} \big| 
 &  \leq C \omega^{-\frac{q}{2}}    \left \| \K^{q/2} \right\|_{L^{1}_1(B_R)}  \\
 &  \leq C \left(\frac{T}{t_\ast}\right)^{- \frac{p}{2}\left(1 - \theta \right) } |B_{R/2}|
\left \| \K^{q/2} \right\|_{\LL^{1}_1(B_R)}^{\theta}. \notag
\end{align}
We next estimate for the third term on the right of~\eqref{e.crabtcha}. Using~\eqref{e.MstrongLp},  we get
\begin{align*}
\lefteqn{ \sum_{j=0}^{k-1} C^{k-1-j} \left( \frac{\sigma}{\beta}\right)^{k-1-j} \left|\left\{ x\in B_{R} : \M_1\!\left( |f|^2\right) (x) > c \sigma \beta^j t_0 \right\} \right|} \qquad  & \\
& \leq Ct_0^{-\frac p2} \sigma^{-\frac p2}\left( \int_{B_R} \left|\M_1\!\left( |f|^2\right) (x) \right|^{\frac p2} \, dx \right) \sum_{j=0}^{k-1} C^{k-1-j}\left( \frac{\sigma}{\beta}\right)^{k-1-j} \beta^{-jp/2}  \\
&  \leq Ct_0^{-\frac p2} \sigma^{-\frac p2} \left \| f^2 \right\|_{L^{p/2}_1(B_R)}^{\frac p{2}}   \sum_{j=0}^{k-1} C^{k-1-j}\left( \frac{\sigma}{\beta}\right)^{k-1-j} \beta^{-jp/2}  \\
& \leq Ct_0^{-\frac p2} \sigma^{-\frac p2} \beta^{-(k-1)p/2} |B_R|
t_\ast^{\frac p{2}},  
\end{align*}
where in the previous line we used 
the fact that, if $c_0$ is sufficiently small, then
\begin{align*}
 \sum_{j=0}^{k-1}C^{k-1-j} \left( \frac{\sigma}{\beta}\right)^{k-1-j} \beta^{-jp/2} 
  &\leq  \sum_{j=0}^{k-1} C^{k-1-j} \left( c_0 \beta^{-\frac{p}{2}}  \right)^{k-1-j} \beta^{-jp/2} \\
  &\leq \beta^{-(k-1)p/2} \sum_{j=0}^{k-1} 2^{j-k+1}  
  \leq 2\beta^{-(k-1)p/2}.
 \end{align*}
To estimate the expression on the last line of the previous string of inequalities, we observe that
\begin{equation*} \label{}
\beta^{p} = c\sigma^{-\frac {2p}{p-2}} = 
c\left(\frac{T}{t_\ast}\right)^{2\nu } 
\left \| \K^{q/2} \right\|_{\LL^{1}_1(B_R)}^{\frac{4\nu}p}
\end{equation*}
and thus, by the definition of $k$,
\begin{equation*} \label{}
t_0^{-\frac p2} \beta^{-(k-1)p/2} \leq \beta^{p} T^{-\frac p2}\leq  c 
T^{- \frac{p}{2} +2\nu } t_\ast^{- 2\nu } 
\left \| \K^{q/2} \right\|_{\LL^{1}_1(B_R)}^{\frac{4\nu}p}.           \\
\end{equation*}
Using this and substituting for $\sigma$ and using the identity
\begin{equation*} \label{}
\theta = \frac p{2m} + \frac{4\nu} p
\end{equation*}
we get
\begin{equation*}  
t_0^{-\frac p2} \sigma^{-\frac p2} \beta^{-(k-1)p/2} |B_R|  t_\ast^{p/2}
 \leq C |B_{R/2}| \left(\frac{T}{t_\ast}\right)^{- \frac{p}{2}\left(1 - \theta \right) }  \left \| \K^{q/2} \right\|_{\LL^{1}_1(B_R)}^{\theta}.
\end{equation*}
We therefore obtain
\begin{multline}
\label{e.boink3}
\sum_{j=0}^{k}C^{k-j} \left( \frac{\sigma}{\beta}\right)^{k-j} \left|\left\{ x\in B_{R} : \M_1\!\left( |f|^2\right) (x) > c \sigma \beta^j t_0 \right\} \right| \\
\leq C |B_{R/2}| \left( \frac{T}{t_\ast}\right)^{- \frac{p}{2}\left(1 - \theta \right) }  
\left \| \K^{q/2} \right\|_{\LL^{1}_1(B_R)}^{\theta}.
\end{multline}
We now insert the inequalities~\eqref{e.boink1},~\eqref{e.boink2} and~\eqref{e.boink3} into~\eqref{e.crabtcha}, taking note that $\beta^kt_0 \leq T$ and $p> 2$, to get 
\begin{align*} 
\frac{\left| B_{R/2}  \setminus A( T) \right|  }{\left| B_{R/2} \right|}
&  \leq C \left(\frac{T}{t_0} \right)^{-\frac p2+\nu}\left \| \K^{q/2} \right\|_{\LL^{1}_1(B_R)}^{\frac{2\nu}p}
  +   C  \left( \frac{T}{t_\ast} \right)^{- \frac{p}{2}\left(1 - \theta \right) }  \left \| \K^{q/2} \right\|_{\LL^{1}_1(B_R)}^{\theta}  \\
& \leq  C  \left( \frac{T}{t_\ast} \right)^{- \frac{p}{2}\left(1 - \theta \right) }  \left \| \K^{q/2} \right\|_{\LL^{1}_1(B_R)}^{\theta} .
\end{align*}
This completes the proof of the proposition.
\end{proof}

\subsection*{Acknowledgements}

The second author was partially supported by Deutsche Forschungsgemeinschaft grant 
no. HO-4697/1-1.

\bibliographystyle{plain}
\bibliography{W1p}

\end{document}